\documentclass[11pt]{amsart}
\usepackage{amscd,amsmath,amssymb,amsthm}
\usepackage{color}
\usepackage{textcomp}
\usepackage{graphicx}
\usepackage{tikz}

\newtheorem{theorem}{{\bfseries Theorem}}[section]

\newtheorem{proposition}[theorem]{{ Proposition}}
\newtheorem{lemma}[theorem]{{ Lemma}}
\newtheorem{corollary}[theorem]{{ Corollary}}
\newtheorem{definition}[theorem]{{ Definition}}
\newtheorem{example}[theorem]{{ Example}}
\newtheorem{remark}[theorem]{{ Remark}}

\newcommand{\bt}{\begin{theorem}}
\newcommand{\et}{\end{theorem}}
\newcommand{\bl}{\begin{lemma}}
\newcommand{\el}{\end{lemma}}
\newcommand{\bp}{\begin{proposition}}
\newcommand{\ep}{\end{proposition}}
\newcommand{\bex}{\begin{example}}
\newcommand{\eex}{\end{example}}
\newcommand{\bc}{\begin{corollary}}
\newcommand{\ec}{\end{corollary}}
\newcommand{\bo}{\begin{proof}}
\newcommand{\eo}{\end{proof}}
\newcommand{\bd}{\begin{definition}}
\newcommand{\ed}{\end{definition}}
\newcommand{\br}{\begin{remark}}
\newcommand{\er}{\end{remark}}
\newcommand{\be}{\begin{enumerate}}
\newcommand{\ee}{\end{enumerate}}

\newcommand{\cD}{\mathcal{D}}

\newcommand{\cR}{{\mathcal{ R}}}

\newcommand{\Z}{{\mathbb Z}}

\newcommand{\N}{{\mathbb N}}
\newcommand{\R}{{\mathbb R}}

\newcommand{\D}{{\mathbf D}}
\newcommand{\EQ}{{\mathbf {EQ}}}
\newcommand{\SR}{{\mathbf {SR}}}
\newcommand{\UR}{{\mathbf {UR}}}

\begin{document}

\title{Strongly Rigid Flows}
\author{Anima Nagar and Manpreet Singh}
\address{Department of Mathematics, Indian Institute of Technology Delhi,
Hauz Khas, New Delhi - 110016, INDIA}

    \vspace{.2cm}
    
    \begin{abstract}
    	We consider flows $(X,T)$,  given by  actions  $(t, x) \to tx$, on a compact metric space $X$ with a discrete $T$ as an acting group. We study a new class of flows - the \textsc{Strongly Rigid} ($ \mathbf {SR} $) \ flows, that are properly contained in the class of distal  ($ \mathbf D $) flows and properly contain the class of all equicontinuous  ($ \mathbf {EQ} $) flows.  Thus, $\mathbf {EQ} \ \text{flows} \subsetneqq \mathbf {SR} \ \text{flows} \subsetneqq \mathbf{ D} \ \text{flows}$.
    
	The concepts of equicontinuity, strong rigidity and distality coincide for the induced flow $(2^X,T)$. We observe that strongly rigid $(X,T)$ gives distinct properties for the induced flow $(2^X,T)$ and  its enveloping semigroup $E(2^X)$. 	We further study strong rigidity in case of particular semiflows $(X,S)$, with $S$ being a discrete acting semigroup.

\end{abstract}

\maketitle

\renewcommand{\thefootnote}{}

\footnote{\emph{keywords:} flows, semiflows, distality, equicontinuity, strongly rigid}


\footnote{The second named author thanks CSIR for financial support.}

 In our study, $(X,T)$ is a flow where $(X,d)$ is  a compact, infinite metric  space and $T$ is an infinite, discrete though not necessarily abelian topological group. When $T = \Z$  we  consider the cascades $(X,f)$, for self-homeomorphism $f$ on $X$.
 
 Identify $t\in T$ with the map $t\longrightarrow tx$. The \emph{enveloping semigroup}  $E(X,T) = E(X)$ of a flow $(X, T)$
 is defined as the closure of the set $\{t : X \to X : t \in T\}$ considered as a subset of $X^X$ given the product
 topology. Equivalently, for the cascade $(X,f)$,  $E(X,f) = E(X)$ is the closure of the family $\{f^n: n \in \Z\}$ in the product topology on $X^X$. 
 
 The minimal idempotents in $E(X)$ play an important role in determining the dynamical properties of $(X,T)$. In this article, we observe that convergence at these minimal idempotents define a new class of flows - the \textsc{Strongly Rigid  Flows}, that are properly contained in the class of distal flows and properly contain the class of all equicontinuous flows. 
 
\bigskip
 
 We  take the following notation for the classes of flows,

 $$ \begin{array}{|c|c|}
 	
 	\hline
 	
 	\text{Equicontinuous Flows}  & \EQ \ \text{flows}\\ 	
 	\hline 
 	\text{Strongly Rigid Flows}  & \SR \ \text{flows} \\ 	
 	\hline 
 	\text{Distal Flows}  & \D \ \text{flows} \\
 	
 	\hline
 	
 \end{array}$$

 \bigskip
  
  and observe the following containment diagram: 
  
  
  \begin{center}
  	\begin{tikzpicture}
  		\draw [thick] (0,-1) arc (-90:270:3cm and 2.2cm);
  		\draw [thick] (0,-0.5) arc (-90:270:2.5cm and 1.5cm);
  		\draw [thick] (0,0) arc (-90:270:1.5cm and .75cm);
  		\node [yshift=.7cm] (0,0) {$\EQ \ \text{flows}$};
  		\node [yshift=2cm] (0,0) {$ \SR \ \text{flows}$};
  		\node [yshift=3cm] (0,0) {$ \D \ \text{flows}$};
  	\end{tikzpicture}
  \end{center}
 
\bigskip

The flow $(X,T)$ induces the flow $(2^X,T)$. For the induced flow $(2^X,T)$, the concepts $\mathbf {EQ} \Leftrightarrow \mathbf {SR} \Leftrightarrow \mathbf D$. Strongly rigid $(X,T)$ gives distinct properties for the induced flow $(2^X,T)$ and  its enveloping semigroup $E(2^X)$.

\bigskip

In Section 1, we recall some basic definitions and results that we shall use. Section 2 is devoted to the study of strongly rigid flows and the dynamical consequences of this property. In Section 3, we extend our study to the case of  semiflows.

\section{Basic Definitions and Results}

Most of the standard definitions and results can be seen in \cite{AUS, ELL, GH}. For more recent results on induced dynamics we follow \cite{AAA, TDES, AN}.

 Consider the flow $(X,T)$. The set $Tx = \{tx: t \in T\}$  is called the \emph{orbit} of a point $x \in X$. A point $x \in X$ is called \emph{recurrent}  if it is an accumulation point of its orbit. A subset $A\subseteq X$ is \emph{invariant} if $TA=\{ta\ :\ a\in A, \ t\in T\}\subseteq A$. If $A$ is \emph{invariant} then $ T\times A \rightarrow A$ is an action and $(A, T)$ is a flow, and is called a \emph{subsystem} of $(X,T)$. 

For flows $(X,T)$ and $(Y,T)$ and $\phi:X \to Y$,  $(Y,T)$ is said to be a \emph{factor} of $(X,T)$ if $\phi$ is a continuous surjection and $\phi(tx)= t\phi(x) \ \forall\ x\in X, t\in T$. These flows are \emph{conjugate} if such a $\phi$ is a homeomorphism.

\bigskip

$(X,T)$ is \emph{point transitive} if there exists a point $x_0\in X$ with $\overline{Tx_0}= X$, where $\overline{A}$ denotes the closure of $A$, and   is called  \emph{(topologically) transitive} if for every nonempty, open subsets $U$ and $V$ in $X$, there exist a $t\in T$ such that $tU\cap V\neq\emptyset$. These two concepts are equivalent on a perfect, compact metric flow.

\bigskip

 A subset $M\subseteq X$ is said to be \emph{minimal} if $M \neq \emptyset $, $M$ is closed, $M$ is invariant and $ M$ is minimal with respect to these properties and the flow $(X,T)$ is \emph{minimal} if $X$ itself is minimal. Equivalently, $M$ is minimal if $\overline{Tx}=M$ for all $x\in M$. A point $x\in X$ is called \emph{almost periodic} if $\overline{Tx}$ is minimal set.  A flow $(X,T)$ is called \emph{pointwise almost periodic} if every $x\in X$ is almost periodic, and is called called \emph{densely almost periodic} if there is a dense set $A \subset  X$ of  almost periodic points. 
 
 \bigskip
 
 The points $x,y \in X$ are \emph{proximal} if the orbit closure of the point $(x, y)$ in the flow $(X \times X, T)$  intersects the diagonal $\Delta$, and the flow $(X,T)$ is proximal if all distinct points are pairwise proximal. The points $x,y \in X$ are \emph{distal} if they are not proximal, and the flow $(X,T)$ is distal if all distinct points are pairwise distal. Distal flows are pointwise almost periodic. $(X,T)$ is said to be \emph{equicontinuous at a point $y\in X$} if for every $\epsilon> 0$ there exists a neighborhood $U$ of $y$ such that for every $x\in U$ and every $t\in T$; we have $d(tx,ty)<\epsilon$ and the flow is \emph{equicontinuous} if it is equicontinuous for every point of $X$. All equicontinuous flows are distal. 

\bigskip

Note that in case of cascades $(X,f)$, one usually talks of forward orbits.  Hence for cascades,  the \emph{orbit} of $x$ is defined as $\mathcal{O}(x) \ = \ \{f^n(x) : n \in \N \}$, and a point $x \in X$ is said to be recurrent if there exists a sequence $\{n_k\}$ with $n_k \nearrow \infty$ such that $f^{n_k}(x) \to x$.  $(X,f)$ is called \emph{topologically transitive}  if for every nonempty, open pair
$U, V \subset X$, there exists a $n \in \N$ such that the set  $f^n(U) \cap V$ is nonempty, with it being \emph{minimal} if $\mathcal{O}(x)$ is dense in $X$ for every $x \in X$.

 A cascade $(Y,g)$ is called a \emph{factor} of the cascade $(X,f)$ if there exists a continuous surjection $h:X \to Y$ such that $h\circ f=g\circ h$. The systems are called \emph{conjugate} if the maps $\phi$ is a homeomorphism. A cascade $(X,f)$ is equicontinuous if the family $\lbrace f^n: n\in \Z \rbrace$ is equicontinuous on $X$.

\bigskip

We follow some definitions given for cascades in \cite{RIG}. A cascade $(X,f)$ is called \emph{weakly rigid} if for any $x_1,x_{2}, \ldots, x_n \in X$ and $\epsilon > 0$, there is a $k \in \N$ such that $d(x_i, f^k (x_i)) < \epsilon,$ $i = 1, \ldots, n$; and is said to be \emph{uniformly rigid} if there is a sequence $\lbrace n_k \rbrace$ with $n_k \nearrow \infty$ such that $f^{n_k} \to $ the identity map uniformly on $X$. 

\bigskip

For any topological space $X$,  if $x\in X$ and $U\subset X$  is open then the set $(x,U)=\lbrace f:X\rightarrow X: f(x)\in U  \rbrace$ forms a sub-basic open set of the \emph{point open topology} i.e. \emph{product topology} on $X^X$, giving pointwise convergence on $X$. If $C\subset X$ is a compact subset, and $U\subset X$ is an open subset,  the set $(C,U)=\lbrace f:X\rightarrow X: f(C)\subset U  \rbrace$ forms a sub-basic open set in the \emph{compact-open  topology} on $X^X$. Compact-open topology gives the uniform convergence on the compact metric space $X$.

On an equicontinuous family $\mathcal{F}$, the compact-open topology is equivalent to the point-open topology.

\bigskip

$2^X$ is  the space of all nonempty closed subsets of $X$,  endowed with the Hausdorff topology.  

\smallskip

Given a point $x \in X$ and a closed set $A \subseteq X$, recall
$d(x, A) = \inf \limits_{a \in A} d(x, a)$ and  \emph{the Hausdorff metric} is defined as 

$$d_H(A,B) \quad = \quad \max \{\sup \limits_{a \in A} d(a, B),  \sup \limits_{b \in B} d(b, A)\}, \ \forall \ A, B \in 2^X.$$

Since $X$ is compact, we occasionally use an equivalent topology on $2^X$. Define for
any collection $\{ U_i : 1 \leq i \leq n \}$ of  open and nonempty subsets of $X$, $\langle U_1, U_2, \ldots U_n \rangle =  \{E \in 2^X :E \subseteq \bigcup \limits_{i=1}^n
U_{i}, \ E \bigcap U_{i} \neq \phi,  \textrm{ } 1 \leq i \leq n \}$.

The topology on $2^X$, generated by such collection as basis, is
known as the {\emph{Vietoris topology}}. 

\bigskip

Note that $(2^X,d_H)$ is also a compact metric space. The induced acting topological group $T$ is defined as: $\forall \ t \in T$, $tA = \{ta:a \in A\}$. Thus we have the induced flow $(2^X,T)$. For cascades $(X,f)$ we denote the induced flow by cascades $(2^X,f_*)$ as in \cite{AAA}.

%
%
%

\bigskip

Let $e \in T$ be the \emph{identity} mapping on $X$. $E(X,T) = E(X)$ is the closure of $T \subset X^X$ given the point-open topology. $E(X)$ forms a semigroup called the \emph{enveloping semigroup}. A non empty subset $I\subset E(X)$ is called a \emph{left ideal} if $E(X)\cdot I\subset I$; i.e. $\alpha \in I$, $p\in E(X) \Rightarrow p \alpha \in I$.  $I$ is a \emph{minimal ideal} if and only if $I$ is closed in $E(X)$ and  does not contain any left ideal as a proper subset. Then an idempotent $u \in I$ is called a minimal idempotent. It is known that $x \in X$ is distal if and only if $ux = x$ for a minimal idempotent $u \in E(X)$.

 \bp \cite{TDES} \label{eq=d} $(X,T)$ is equicontinuous if and only if  $(2^{X},T)$  is equicontinuous if and only if $(2^{X},T)$  is distal. \ep

We give a proof here for the sake of completion.
 
 \vskip .5cm
 
 \bo  It is simple to see that $(X,T)$ being equicontinuous implies that $(2^X, T)$ is also equicontinuous and that further implies that $(2^X,T)$ is distal. We will show that $(2^X,T)$ if distal implies that $(X,T)$ is equicontinuous.
 
 For this we consider the  \emph{regionally proximal relation} $RP \subset X \times X$ for $(X,T)$.
 We say that $(x,y) \in RP$ if
 there are nets $\{x_n\}$ and $\{y_n\}$ in $X$ with $x_n \to x$ and $y_n \to y$, $z \in X$
 and net $\{t_n\}$ in $T$ such that $({t_n}(x_n),{t_n}(y_n))\to (z,z)$ in $X \times X$.  It is
 known  that $(X,T)$ is equicontinuous if and only if $RP$ is the identity relation in $X \times X$.
 
 Suppose  $(X,T)$ is not equicontinuous.
 Then there are $x \neq y \in X$ with  $(x,y) \in RP$. Then
 we have nets $x_n \to x$, $y_n \to y$ and $t_n \in T$ with
 $({t_n}(x_n),{t_n}(y_n))\to (z,z)$ for some $z \in X$.
 Let $C=\{x_1,x_2, \dots\} \cup \{x\}$ and $D=\{y_1,y_2, \dots\} \cup \{y\}$. We
 can always take $C \cap D=\emptyset$. Since  $(2^X,T)$ is distal, the orbit closure
 of $(C,D)$ in $2^X \times 2^X$ should be minimal. Let $(C',D')$ be the limit point of the set $\{ ({t_n}(C ,D) \}$ in $2^X \times 2^X$. But $z \in C' \cap D'$, and so we have $C' \cap D' \neq \emptyset$. Now if the $T- $ orbit closure of $(C,D)$ were minimal, then $(C,D)$ would be in
 the $T -$ orbit closure of $(C',D')$. But that would imply that  $C \cap D \neq \emptyset$, a contradiction. \eo

We recall from \cite{AN}; for $p \in E(X)$ and $A \in 2^X$,
$$\cD_p(A) = \lbrace y \in X: \exists\ \text{nets} \ \lbrace a_{i}\rbrace \ \text{in} \ A, \ \text{and} \ \lbrace t_{i} \rbrace \ \text{in} \ T \ \text{with} \ t_i \to p  \ \text{such that}  \ t_{i}a_{i}\rightarrow y\rbrace.$$

 $\cD_p: 2^X \to 2^X$ is a function, and for each   $p \in E(X)$ and a finite $A \in 2^X$, $\cD_p(A) = pA$. For $t \in T$ and $A \in 2^X$, $\cD_t(A) = tA$.  Also  $\cD_u$ is an idempotent function for every idempotent $u \in E(X)$.

Note that $A \in 2^X$ is an almost periodic point for $(2^X, T)$ if  $\cD_u(A) = A$ for a minimal idempotent $u   \in E(X)$. The almost periodic points for $(2^X, T)$ are precisely the elements of $2^X$ in the range of $\cD_u$, \ for every minimal idempotent $u  \in E(X)$. 

\bigskip

Lastly, recall that the convergence in $E(2^X)$ is the uniform convergence on $X^X$.

\bigskip

\begin{section}{ Flows}
	
In $E(X)$ any convergence  is pointwise, and every minimal idempotent is an accumulation point for $E(X)$. We can consider the possibility of the minimal idempotents to be accumualtion points for $E(X)$ even when endowed with the stronger uniform topology. This defines a new class of flows:

\bd
A flow $(X,T)$ is said to be \emph{strongly rigid} ($ \SR $) if  every minimal idempotent $u \in E(X)$ is an accumulation point in  $E(X)$ when endowed with  the compact-open topology, i.e. there exists a (non-trivial) net $\lbrace t_{k} \rbrace$ in $T$ such that the convergence $t_{k} \longrightarrow u$  is uniform.
\ed

\begin{remark}
All $ \EQ $ flows are $ \SR $.
\end{remark}

\bp \label{d} Every $ \SR $ flow is $ \D $. \ep
\begin{proof}
If $(X,T)$ is strongly rigid then for any minimal idempotent $u\in E(X)$, there is a (non-trivial) net $\lbrace t_{k} \rbrace$ in $T$ with $t_{k}\longrightarrow u$ and the convergence is uniform. As $X$ is compact, we observe that $t_{k}(X)\longrightarrow u(X)$, and since $t_{k}(X)=X $ we have $u(X)=X$. So $u(x)=x$ for all $x\in X$ and this means that $u=e$(identity) on $X$. Thus $e$ is the only minimal idempotent in $E(X)$ and so $(X,T)$ is distal. 
\end{proof}

\br We review the definition of uniformly rigid   given in \cite{RIG}. A flow $(X,T)$ is said to be \emph{uniformly rigid} ($ \UR $) if the identity  $e \in E(X)$ is an accumulation point in  $E(X)$ when endowed with  the compact-open topology, i.e. there exists a (non-trivial) net $\lbrace t_{k} \rbrace$ in $T$ such that the convergence $t_{k} \longrightarrow e$  is uniform.

Note that $\UR$ flows need not be $\SR$  since $e$ need not always be a minimal idempotent. However for flows, $ \SR = \D \cap \UR$. \er

The converse of Proposition \ref{d} is not true. We have a $ \D $ cascade which is not $ \SR $. We borrow the below example from \cite{RIG}:

\bex \label{A1}
Let $X=\lbrace re^{i\theta}: \ r=1 - \frac{1}{2^n}, \ n \in \N, \ \text{and} \ r=1, \ 0\leq\theta\leq 2\pi \rbrace$ and $f:X\rightarrow X, f(re^{i\theta})=re^{\theta+2\pi ir}$. The cascade $(X,f)$ is distal but not uniformly rigid. Hence is not strongly rigid.
\eex

Also there exists a $ \SR $ cascade which is not $ \EQ $. We borrow the below example from \cite{OP}, and discuss some details here for the sake of completion:

\bex \label{A2}
$X= \R/\Z \times \R/\Z$. For   $n_{1}=100, \ n_{k+1}={(n_{1}n_{2}\ldots n_{k})}^{3}$, \ $k \geq 1$;  take the irrational number $\alpha=\sum\limits_{j=1}^{\infty}\frac{1}{n_{j}}$.  

Let $\phi:\R/\Z \rightarrow \R/\Z$ be defined as $\phi(x)=\sum\limits_{k=-\infty}^{+\infty}(e^{2 \pi i n_{k}\alpha}-1)e^{2 \pi i n_{k}x}$.

For every $x \in \R/\Z$, it can be seen that $\phi(x)\in \R/\Z$ and  the function $\phi:\R/\Z \rightarrow \R/\Z$  is continuous.

 The map $T:X\rightarrow X$ is defined as;

\begin{center}
	$T(x,y)=(x+\alpha, y+\phi(x)) (\mod1).$
\end{center}

Note that $T$ is surjective on $X$.

Then as shown in \cite{OP}, it can be seen that the cascade $(X,T)$ is distal. 

For an $n\in \Z$, $T^{n}(x,y)=(x+n\alpha, \sum\limits_{k=1}^{n-1} \phi(x+k\alpha))(\mod 1)$

Suppose  $(x_{1},y_{1})$, $(x_{2},y_{2}) \in X$ is a proximal pair. Then there is a sequence $\lbrace m_{k} \rbrace$ in $\Z$ such that 
\begin{center}
	$T^{m_{k}}(x_{1},y_{1})\longrightarrow (a,b) \longleftarrow T^{m_{k}}(x_{2},y_{2})$ 
\end{center}

Then, 
	$ x_{1}+m_k\alpha(\mod 1)  \longrightarrow a \longleftarrow x_{2}+m_k \alpha (\mod 1)$ and 
	
	$y_{1}+ \sum\limits_{k=1}^{m_k-1} \phi(x_{1}+k\alpha)\longrightarrow b \longleftarrow y_{2}+ \sum\limits_{k=1}^{m_k-1} \phi(x_{2}+k\alpha) (\mod 1)$.

So 	 $ x_{1}=x_{2}$ and $ y_{1}=y_{2}$, giving $(x_{1},y_{1})=(x_{2},y_{2})$.

Also it is shown in \cite{OP} that $(X,T)$ is not equicontinuous. For $(x_{n_{l}},0)\longrightarrow (0,0),(l\longrightarrow \infty)$ with $x_{n_{l}}=\frac{1}{n_{1}n_{l}}$, and  $m_{i}=\frac{n_{l}^{3}}{n_{1}}$, $\delta=\frac{1}{1000}$, it is shown there that 
	$d(T^{m_{l}}(x_{n_{l}},0),T^{m_{l}}(0,0))>\delta$.

And  as shown in \cite{OP}, for every $(x,y) \in X$, $T^{n_s}(x,y) \to (x,y)$ uniformly.

Thus the cascade $(X,T)$ is distal, non-equicontinuous and strongly rigid. 

Hence $(X,T)$ is $ \SR $ and not $ \EQ $.

\eex

Can we have a more stronger version of the concept of $ \mathbf{SR} $? The definition does give such a possibility:

\bd
A $ \mathbf{SR} $ flow is called \emph{very strongly rigid}$\ (\mathbf{VSR}) $ if every net $\lbrace t_{k} \rbrace$ in $T$ such that $t_{k}\longrightarrow e$ in $E(X)$ implies that $t_{k}\longrightarrow e$ uniformly on $X$. 
\ed

This means that very strongly rigid flows are those strongly rigid flows where convergence to identity is always uniform. But that leads to the below observation:

\bp
$(X,T)$ is  $ \mathbf{VSR} $ if and only if $(X,T)$ is  $ \EQ $.
\ep
\begin{proof}
	If $(X,T)$ is equicontinuous then trivially it is very strongly rigid.
	
	Conversely if $(X,T)$ is very strongly rigid then since every $t_{k}\longrightarrow e$ is uniform, and so $t_{k}A\longrightarrow A$ for all $A\in 2^X$, hence $$D_{e}(A)  = \lbrace y \in X: \exists\ \text{nets} \ \lbrace a_{i}\rbrace \ \text{in} \ A, \ \text{and} \ \lbrace t_{i} \rbrace \ \text{in} \ T \ \text{with} \ t_i \to e  \ \text{such that}  \ t_{i}a_{i}\rightarrow y\rbrace =A$$ and so all $A\in 2^X$ are almost periodic i.e $(2^X,T)$ is distal and hence equicontinuous. This is equivalent to $(X,T)$ being equicontinuous by Proposition \ref{eq=d}.   
\end{proof}

So we have the following strict inclusions, with no intermediate classes: 
\begin{center}
 $ \EQ \ \text{flows}  \subsetneqq \SR \ \text{flows} \subsetneqq \D \ \text{flows}$.
\end{center}

Though we find that this relation can be extended. We define analogous to \cite{RIG}:

\bd A flow $(X,T)$ is called \emph{weakly rigid} ($ \mathbf{WR} $) if for any $x_1,x_{2}, \ldots, x_n \in X$ and $\epsilon > 0$, there is a $t \neq e \in T$ such that $d(x_i, t x_i) < \epsilon,$ $i = 1, \ldots, n$. \ed

As observed in \cite{RIG, TDES}, we see that $ \mathbf{WR} $ flows are precisely those flows that have a  non-isolated identity element in their enveloping semigroup. Thus, we can say:

\bp \label{wr} A $\D$ flow is always $ \mathbf{WR} $. \ep

The converse to Proposition \ref{wr} does not hold true and for that we have a counter-example in \cite{RIG}. The class of $ \mathbf{WR} $ flows is very large and we will not investigate it.

\vspace{0.5cm}

\subsection{Some Dynamical Properties}

\bp \label{sr-factor}
 A factor of $\SR$ flow is $\SR$.
\ep
\begin{proof} Let $(Y,T)$ be a factor of the flow $(X,T)$ and $(X,T)$ be strongly rigid. 	 Let $\phi:(X,T)\rightarrow (Y,T)$ be the factor map. Then there is a factor map $\Phi:E(X)\rightarrow E(Y)$ defined as $\Phi (p)(\phi(x))= \phi(px)$. Hence $(Y,T)$ is distal because $(X,T)$ is distal.
	
	Now since  $(X,T)$ is strongly rigid, there is a (non-trivial) net $\lbrace t_{i} \rbrace$ in $T$ such that $t_{i}\longrightarrow e$ uniformly. Since $\Phi$ is a factor map, $\Phi(t_{i})\longrightarrow \Phi(e)=e$ the identity in $E(Y)$. 
	
	 Since $\Phi$ is uniformly continuous with $X$ being compact,  for a given $\epsilon>0$, there is a $\delta>0$ such that $d(\phi(x),\phi(y))<\epsilon$ where $d(x,y)<\delta$. Now for this $\delta$, there is an indexing set $\Gamma$  such that for all $x\in X$, $d(t_{i}(x),x)<\delta$ for all $i \in \Gamma$, which gives $d(\phi(t_{i}(x),\phi(x))<\epsilon$ for all $i \in \Gamma$ and for all $x\in X$. Thus $ d(\Phi(t_{i})(\phi(x)),\phi(x))<\epsilon$.  
	
	Let $\phi(x)=y \in Y$ and so $d(\Phi(t_{i})(y),y)< \epsilon$ for all $i \in \Gamma$ and $y \in Y$ as $\phi$ is surjective. Thus, $\Phi(t_{i})\longrightarrow \Phi(e)=e$ uniformly and so $(Y,T)$ is strongly rigid.   
\end{proof}

Suppose $(X,T)$ is a flow and $Y\subset X$ is closed and invariant set, then we note that $E(X,T)\subset E(Y,T)$. 

\bp
A subsystem of an $\SR$ flow is also $\SR$. 
\ep
\begin{proof}
	Suppose $(X,T)$ strongly rigid, and $Y\subset X$ is closed and invariant.  $(X,T)$ is distal and so $(Y,T)$ is also distal. So $E(X,T)$ is subgroup of the group $E(Y,T)$. There is a net $\lbrace t_{k} \rbrace$ in $T$ such that $t_{k}\longrightarrow e$ uniformly on $X$. Since $e$ is the only idempotent in $E(Y)$, $(Y,T)$ is also strongly rigid.
\end{proof}

Recall that for an indexing set $\Gamma$ and flows $(X_{\gamma},T)$  for $\gamma \in \Gamma$, $ (\prod\limits_{\gamma \in \Gamma} X_{\gamma},T)$ is distal if and only if $(X_{\gamma},T)$ is distal for each $\gamma \in \Gamma$.

\bp \label{product}
 $(X,T)$ is $\SR$ if and only if for any integer $k>0$, $(X^k,T)$ is $\SR$.
\ep
\begin{proof}
	Since $(X,T)$ is factor of $(X^k, T)$ under the projection map. So by Proposition \ref{sr-factor}  $(X, T)$ is strongly rigid whenever for any integer $k>0$, $(X^k,T)$ is. 
	
	Conversely, if $(X,T)$ is strongly rigid, then $(X,T)$ is distal. Therefore $(X^k,T)$ is also distal for any $k > 0$. Since $E(X^k)\cong \Delta E(X)^k$, so the identity $(e)=(e,e,\ldots, e)\in E(X^k)$ is the only idempotent element here. Since $(X,T)$ is strongly rigid, there is a (non-trivial) net $\lbrace t_{j} \rbrace$ in $T$ such that $t_{j} \longrightarrow e$ uniformly, which gives that $(t_{j})\longrightarrow (e)$ uniformly in $E(X^k)$ and hence $(X^k,T)$ is strongly rigid.    
\end{proof}

\bt \label{subnet}
$(X,T)$ is $ \SR $ if and only if the net $\lbrace t_{k} \rbrace$ in T with $t_{k}\longrightarrow e$ uniformly has a subnet converging to $e$ indexed by the directed set $\mathcal{N}_{e}$, the neighbourhood basis of $e$ in the compact-open topology.
\et
\begin{proof}
	Suppose $(X,T)$ is strongly rigid, so for an indexing set $\Lambda$ there is a  $P:\Lambda \rightarrow T$ with $P(k)=t_{k}$ such that $t_{k}  \longrightarrow e$ uniformly. Now for any neighbourhood $N$ of identity $e$ in the compact-open topology on $X^X$, there is an $r \in \Lambda$ such that $t_{k}\in N$ for all $k\geq r$ and for $k<r$, $t_{k} \notin N$. Let $s= \inf\lbrace r\in \Lambda: t_{r}\in N\rbrace$. Such an infimum always exists since $\Lambda$ is an indexing set for a net converging to $e$ uniformly. So we define a function $\phi:\mathcal{N}_{e}\rightarrow \Lambda$  as $\phi(N)=s= inf\lbrace r\in \Lambda: t_{r}\in N\rbrace$. This map is well defined.
	\vspace{1mm}
	
	\underline{\bf{claim 1}: $\phi$ is increasing}: Let $N_{m} \subset N_{l}$ in $\mathcal{N}_{e}$, then $N_{l} \leq N_{m}$. So for $k\geq \phi(N_{l})=l^{\prime}$, $t_{k}\in N_{l}$ and for $k\geq \phi(N_{m})=m^{\prime}$, $t_{k}\in N_{m}$. Since $N_{m}\subset N_{l}$, $inf\lbrace s\in \Lambda: t_{s}\in N_{l}\rbrace \leq inf\lbrace s\in \Lambda: t_{s}\in N_{m}\rbrace\Rightarrow \phi(N_{l})\leq \phi(N_{m})$.
	
	\vspace{1mm}
	\underline{\bf{claim 2}: $\phi$ is cofinal}: Let $\lambda \in \Lambda$. For an $N_{p} \in \mathcal{N}_{e}$ let $\phi(N_{p})=p^{\prime}$. If $p^{\prime}<\lambda$,  consider the open neighbourhood $N_p \setminus \lbrace t_{k_{0}}: k_{0}> \lambda\rbrace = {N'_{p}}$(say). Then $\inf \lbrace s\in \Lambda: t_{s}\in {N'_{p}}\rbrace > \lambda$. Therefore $\phi({N'_{p}})= \inf  \lbrace s\in \Lambda: t_{s}\in {N'_{p}}\rbrace > \lambda$ and so $\phi$ is cofinal.
	
	\vspace{1mm}
	 Hence $\phi \circ P:\mathcal{N}_{e}\rightarrow X$ gives a subnet of $\lbrace t_{k} \rbrace$ with $t_{\phi(N_{l})} \longrightarrow e$. 
	
	\vspace{2mm}
	The converse follows from  definition.
\end{proof}

\br We can take an alternate definition of $\SR$ flows as those distal flows for which there is a net converging to $e$ indexed by the directed set $\mathcal{N}_{e}$, the neighbourhood basis of $e$ in the compact-open topology. \er

This leads to:

\bp
$(X,T)$ and $(Y,T)$ are $\SR$ flows if and only if $(X\times Y, T)$ is $\SR$. 
\ep
\begin{proof}  If  $(X\times Y, T)$ is strongly rigid then by Proposition \ref{sr-factor} both $(X,T)$ and $(Y,T)$ are strongly rigid.
	
	\vspace{2mm}
	
	Conversely, suppose $(X,T)$ and $(Y,T)$ are strongly rigid flows then there are nets $\lbrace t_{k} \rbrace$ and $\lbrace t_{\alpha} \rbrace$  such that $t_{k}\longrightarrow e$ and $t_{\alpha}\longrightarrow e$  in $X^X$ and $Y^Y$ respectively, given the compact open topology. By Theorem \ref{subnet} these nets have converging subnets with the directed set as neighbourhood basis of $e$. Suppose $\phi_{1}:\mathcal{N}_{e}\rightarrow T$ and $\phi_{2}:\mathcal{N}_{e}\rightarrow T$ give these subnets converging to $e$. So for any $\epsilon>0$, there are $N_{1},N_{2}\in \mathcal{N}_{e}$ such that		$d(t_{\phi_{1}(N_{l})}(x),x)<\epsilon$ for all $N_{l}\geq N_{1}$,   $x\in X$, and $d(t_{\phi_{2}(N_{\gamma})}(y),y)<\epsilon$ for all $N_{\gamma}\geq N_{2}$, \ $y\in Y$. Consider $$\mathcal{S}_{e}=\lbrace N_{l}\cap N_{\gamma}: N_{l},N_{\gamma}\in \mathcal{N}_{e}\ \text{such that}\ \phi_{1}(N_{l}) \geq \phi_{1}(N_1), \phi_{2}(N_{\gamma}) \geq \phi_{2}(N_2) \rbrace \subset \mathcal{N}_{e}$$ as a  directed set and the subnet $\phi:\mathcal{S}_{e}\rightarrow T$ defined on the product $(X\times Y)^{X\times Y}$ as $t_{\phi(N_{l}\cap N_{\gamma})}(x,y)=(t_{\phi_{1}(N_{l}\cap N_{\gamma})}x, t_{\phi_{2}(N_{l}\cap N_{\gamma})}y)$.  
	
	Since in the directed set $\mathcal{N}_{e}$, if $N_{i}\leq N_{j}$ then $N_{j}\subset N_{i}$. So $N_{1}\cap N_{2}\subset N_{1}$ and $N_{1}\cap N_{2}\subset N_{2}$ which means $N_{1}\cap N_{2}\geq N_{1},N_{2}$. So for any $(x,y)\in X\times Y$ and $N_{l}, N_{\gamma} \geq N_{1}, N_{2}$, $N_{1}\cap N_{\gamma}\geq N_{1}\cap N_{2}$. For the product metric $D$ on $X\times Y$,
	
	\begin{align*}
		D(t_{\phi(N_{l}\cap N_{\gamma})}(x,y),(x,y)) &=D((t_{\phi_{1}(N_{l})}(x),t_{\phi_{2}(N_{\gamma})}(y)),(x,y))\\ & = \max \lbrace d(t_{\phi_{1}(N_{l})}(x),x), d(t_{\phi_{2}(N_{\gamma})}(y),y) \rbrace \\
		& < \epsilon
	\end{align*}
	
	So the subnet given by $\phi$ converges uniformly to $e$. Hence $(X\times Y,T)$ is strongly rigid.
\end{proof}

\bc For $ \SR $ cascades $(X,f)$ and $(Y,g)$, there exists   sequences $\lbrace n_{k} \rbrace$ with $n_k \nearrow \infty$ and $\lbrace m_{k} \rbrace$ with $m_k \nearrow \infty$ such that $f^{n_{k}} \times g^{m_{k}} \to f^0 \times g^0$ uniformly on $X \times Y$, implying that $(X \times Y, f \times g)$ is $\SR$. \ec

\bigskip

Next observe that $\SR$ flow $(X,T)$ is $\D$ and so is pointwise almost periodic. Thus $X$
can be realized as a union of minimal orbits. 
\medskip

\bd For a flow $(X,T)$,  the \emph{rigidity relation} $\cR \subset X \times X$ is defined as;

\begin{center}
	$\cR= \left\{\begin{array}{rcl}
		& \mbox{there exist a net} \  \lbrace t_{k} \rbrace \ \mbox{with} \ t_{k}\longrightarrow e  \ \mbox{uniformly}, \\
		(x,y): & \mbox{together with subnet} \ \lbrace t_{k}^{\prime} \rbrace\ \mbox{and nets}\  x_{k}\longrightarrow x,\ y_{k} \longrightarrow y\\ & \mbox{such that}\ d(t_{k}^{\prime}x_{k},t_{k}^{\prime}y_{k})\longrightarrow 0  
	\end{array}\right\}.$
\end{center} \ed

\bigskip

\bp For a $\SR$ flow $(X,T)$,  the {rigidity relation} $\cR = \Delta$, the identity relation in $X \times X$. \ep

\bo For any pair $(x,x)\in \Delta$, since $(X,T)$ is strongly rigid, $t_{k}x\longrightarrow x$. So let $x_{k}=t_{k}x$ then for the subnet $\{t_{k}\}$ itself, $d(t_{k}x,t_{k}x)\longrightarrow 0$. So $\Delta \subset \cR$.

Conversely, if $(x,x^{\prime})\in \cR$ then there is $x_{k}\longrightarrow x$ and $x_{k}^{\prime}\longrightarrow x^{\prime}$ and  a net \  $\lbrace t_{k} \rbrace $ \ with \ $ t_{k}\longrightarrow e $  \ {uniformly},\ together with subnet  \ $\lbrace t_{k}^{\prime} \rbrace$ such that $d(t_{k}^{\prime}x_{k},t_{k}^{\prime}x_{k}^{\prime})\longrightarrow 0$. Since $t_{k}\longrightarrow e$ uniformly,  for all $x\in X$,  $d(t_{k}x,x) \longrightarrow 0$. Since $x_{k}\longrightarrow x$ and $x_{k}^{\prime}\longrightarrow x^{\prime}$, eventually  $t_{k}^{\prime}x_{k}\longrightarrow t_{k}^{\prime}x$ and $t_{k}^{\prime}x_{k}^{\prime}\longrightarrow t_{k}^{\prime}x^{\prime}$. So
	$d(t_{k}^{\prime}x, t_{k}^{\prime}x_{k})\longrightarrow 0$, $d(t_{k}^{\prime}x_{k}^{\prime},t_{k}^{\prime}x^{\prime})\longrightarrow 0$. 
	
	Therefore, 
	$d(t_{k}^{\prime}x,t_{k}^{\prime}x^{\prime})\longrightarrow 0$, and

\begin{center}
	$d(x,x^{\prime})\leq d(x,t_{k}^{\prime}x)+d(t_{k}^{\prime}x,t_{k}^{\prime}x^{\prime})+ d(t_{k}^{\prime}x^{\prime},x^{\prime}) \longrightarrow 0$   
\end{center}
Hence $x=x^{\prime}$ and $\cR= \Delta$. \eo

\br We recall the {regionally proximal relation} $RP \subset X \times X$ for $(X,T)$.
Note that $(x,y) \in RP$ if
there are nets $\{x_n\}$ and $\{y_n\}$ in $X$ with $x_n \to x$ and $y_n \to y$, and net $\{t_n\}$ in $T$ such that $d({t_n}(x_n),{t_n}(y_n)) \longrightarrow 0$.  It is
known  that $(X,T)$ is equicontinuous if and only if $RP  = \Delta$.

We note that in general $\cR \subseteq RP$, and for $\EQ$ flows we have $\cR = RP$. \er

\br It is very evident that $\cR = \emptyset$ for $\D$ flows outside the class of $\SR$ flows. \er

\vspace{8mm}

\subsection{Induced Flows}

\vspace{0.5cm}
The above discussions vary in case of induced flows.

Recall Proposition \ref{eq=d}. We  now extend it as:

\bt 
For a flow $(X,T)$,  the following are equivalent;
\begin{enumerate}
	\item $(X,T)$ is $\EQ$.
	\item $(2^X,T)$ is $\EQ$.
	\item $(2^X,T)$ is $\SR$.
	\item $(2^X,T)$ is $\D$.
	\end{enumerate}
\et

Thus we have the following diagram:

$$\begin{matrix}
	(X,T) \  \text{is} \ \EQ & \Leftrightarrow & (2^X,T) \  \text{is} \ \EQ\\
	\Downarrow &      &  \Updownarrow \\
	(X,T) \  \text{is} \ \SR &     &  (2^X,T) \  \text{is} \ \SR \\
	\Downarrow  &     &  \Updownarrow \\
	(X,T) \  \text{is} \ \D  &     &  (2^X,T) \  \text{is} \ \D \\
\end{matrix}$$

\bigskip

For the induced flows though $\EQ$ is equivalent to $\SR$, we have some other observations.

\bp \label{sr-wr} $(X,T)$ is $\SR$ $\Longrightarrow$ $(2^X, T)$ is $\mathbf{WR}$. \ep

The proof follows from Theorem \ref{subnet}  and the fact that the compact-open topology on $X^X$ is equivalent to the point-open topology on ${(2^X)}^{2^X}$.

\br Note that for $\D$ flow $(X,T)$, the induced flow $(2^X,T)$ need not be $\mathbf{WR}$. 

Recall Example \ref{A1}. Let $A=\lbrace re^{2 \pi i}: \ r=1 - \frac{1}{2^n}, \ n \in \N \ \text{and} \ r=1 \rbrace $, and $B=\lbrace re^{2 \pi r i}: \ r=1 - \frac{1}{2^n}, \ n \in \N \ \text{and} \ r=1  \rbrace$ $\subset X$. We notice that for any $0 < \epsilon < \frac{1}{2^{100}} $ there is no $k \in \N$ such that $d_H(A,f_*^k(A)) < \epsilon$ and $d_H(B,f_*^k(B)) < \epsilon$ together.\er

\bigskip

We recall from \cite{AN}, a finite $A \in 2^X$ is an almost periodic point in $(2^X, T)$ if and only if ${uA} = \cD_u(A) = A$ for a minimal idempotent $u \in E(X) $. If $(X, T)$ is a  distal flow, then $A \in 2^X$ is an almost periodic point if $ \cD_e(A) = A$ for the only minimal idempotent $e \in E(X)$. If $(X,T)$ is  distal then each finite set in the induced system
$(2^X,T)$  is  almost periodic and these points are dense in $2^X$. Thus,

\bp For a $\D$ flow $(X,T)$, $(2^X,T)$  is densely almost periodic. \ep

\textit{What are all the almost periodic points in $(2^X,T)$ for a $\D$ flow $(X,T)$?}

We can say something partially here. Let closed $Y \subset X$ be such that the subflow $(Y,T)$ is $\EQ$. Then $(2^Y,T)$ is also $\EQ$. Since $\EQ \Leftrightarrow \mathbf{VSR}$, we note that $\cD_e(A) = A$ for every $A \in 2^Y$. In particular $\cD_e(Y) = Y$.

\bigskip

Again note that for a net $\{t_k\}$ in $T$ with $t_k \to e$ uniformly, every $A \in 2^X$ is an accumulation point for $TA$ in $2^X$. This essentially gives,

\bp For a $\SR$ flow $(X,T)$,  every $A \in 2^X$ is a recurrent point in the induced flow $(2^X,T)$. \ep

\bp  $(X,T)$ is \textbf{SR} $\Longrightarrow$ the rigidity relation $\mathcal{R}=\Delta$ for  $(2^X, T)$. \ep

\bo Since $(X,T)$ is $ \mathbf{SR} $, there is a (non-trivial) net $\lbrace t_{k}\rbrace$ in $T$ such that $t_{k}\longrightarrow e$  in the compact-open topology on $E(X)$.  So $t_{k}\longrightarrow e$ uniformly in $E(2^X)$.

Now for any $A\in 2^X$, $t_{k}A\longrightarrow A$. So for the net $A_{k}=t_{k}(A)$  in $2^X$ and $d_{H}(A_{k},A)\longrightarrow 0$. Hence,  $(A,A)\in \mathcal{R}$. 

Conversely, suppose $(A,B)\in \mathcal{R}$ then there are nets $A_{k}\longrightarrow A$ and $B_{k}\longrightarrow B$ and a subnet $\{t_{k}^{\prime}\}$ of $\{t_k\}$ such that  $d_{H}(t_{k}^{\prime}A_{k},t_{k}^{\prime}B_{k})\longrightarrow 0$. 

Since $A_{k}\longrightarrow A$ and $B_{k}\longrightarrow B\Rightarrow t_{k}^{\prime}A_{k}\longrightarrow t_{k}^{\prime}A$ and $t_{k}^{\prime}B_{k}\longrightarrow t_{k}^{\prime}B$. So
$d_{H}(t_{k}^{\prime}A, t_{k}^{\prime}A_{k})\longrightarrow 0$, and $d_H(t_{k}^{\prime}B_{k},t_{k}^{\prime}B)\longrightarrow 0$. Therefore, 

\begin{center}
	$d_{H}(t_{k}^{\prime}A,t_{k}^{\prime}B)\leq d_H(t_{k}^{\prime}A, t_{k}^{\prime}A_{k})+ d_H(t_{k}^{\prime}A_{k},t_{k}^{\prime}B_{k})+d_H(t_{k}^{\prime}B_{k},t_{k}^{\prime}B) \longrightarrow 0.$
\end{center}

Hence $A=B$ and  $\cR= \Delta$. \eo
.

\bigskip

For  $n\in \N$, let $\it{F}_{n}(X)= \lbrace A\subset X: A\ \text{is closed and card(A)}\leq n \rbrace$ be \emph{the space of finite sets} in $X$. If $(X,T)$ is a flow then $(\it{F}_{n}  
(X),T)$ is also a flow with $tA=\lbrace tx: x\in A \rbrace$. Note that $\it{F}_{n}(X) \subset 2^X$ is also a metric space with the Hausdorff metric $d_H$.

\vspace{0.5cm}
Since $\it{F}_{1}(X)\equiv X$,  we  take $n\geq 2$ always. 

\bp
For a flow $(X,T)$ and $n\in \N$, $(\it{F}_{n}(X),T)$ is a factor of $(X^n,T)$.
\ep
\begin{proof}
Define $\phi:X^n \rightarrow \it{F}_{n}(X)$ as $\phi((x_{1},x_{2},\ldots,x_{n}))=\lbrace x_{1},x_{2},\ldots,x_{n} \rbrace$. Clearly for any $t\in T$, 
\begin{center}
	$\phi(t(x_{1},x_{2},\ldots,x_{n}))= t\phi((x_{1},x_{2},\ldots,x_{n}))$
\end{center}

Let $\lbrace x_{1},x_{2},\ldots,x_{r} \rbrace \in \it{F}_{n}(X)$, $r\leq n$, then $(x_{1},x_{2},\ldots,x_{r},\underbrace{x_{r},x_{r},\ldots, x_{r}}_{n-r\ \text{times}})\in X^n$ and

 $\phi((x_{1},x_{2},\ldots,x_{r},\underbrace{x_{r},x_{r},\ldots, x_{r}}_{n-r\ \text{times}}))= \lbrace x_{1},x_{2},\ldots,x_{r} \rbrace$. So $\phi$ is surjective. 

To check the continuity, let $\lbrace (x^{1}_{i},x^{2}_{i},\ldots, x^{n}_{i}) \rbrace$ be a sequence in $X^n$ such that $(x^{1}_{i},x^{2}_{i},\ldots, x^{n}_{i})\longrightarrow (x_{1},x_{2},\ldots,x_{n}) \Rightarrow x^{j}_{i}\longrightarrow x_{j}$ for each $j=1,2,\ldots,n$ in $X$. So $\lbrace x^{1}_{i},x^{2}_{i},\ldots, x^{n}_{i} \rbrace \longrightarrow \lbrace x_{1},x_{2},\ldots,x_{n} \rbrace$ in $\it{F}_{n}(X)$. Thus,  $\phi((x^{1}_{i},x^{2}_{i},\ldots, x^{n}_{i}))\longrightarrow \phi((x_{1},x_{2},\ldots,x_{n}))$ which means $\phi$ is continuous and hence $(\it{F}_{n}(X),T)$ is a factor of $(X^n,T)$.

\end{proof}

\bc \label{prod_distal}
$(X,T)$ is $\D$ if and only if  $(\it{F}_{n}(X),T)$, $n \in \N$, is $\D$. 
\ec
\begin{proof}
$(X,T)$ is distal $\Leftrightarrow(X^n,T)$ is distal and since $(\it{F}_{n}(X),T)$ is a factor of $(X^n,T)$, $(\it{F}_{n}(X),T)$ is also distal.

Conversely,  assume that $x_{1} \neq x_{2} \in X$ are proximal. But then $\lbrace x_{1} \rbrace$ and $\lbrace x_{2} \rbrace$ are proximal in $(\it{F}_{n}(X),T)$, contradicting that $(\it{F}_{n}(X),T)$ is distal. 

Thus, $x_{1}=x_{2}$ and hence $(X,T)$ is distal.   
\end{proof}

\bc
 $(X,T)$ is $\SR$ $\Leftrightarrow$ $(\it{F}_{n}(X),T), \ n\in \N$ is $\SR$.
\ec

\bc
$(X,T)$ is $\EQ$ $\Leftrightarrow$ $(\it{F}_{n}(X),T), n \in \N$ is $\EQ$.
\ec
	
	In both cases, the proof follows from the fact that each $\it{F}_{n}(X) \subset X^n$, the properties of being $\EQ$ or $\SR$ is  preserved by taking finite products and factors; and the fact that $X \equiv \it{F}_{1}(X) \subset \it{F}_{n}(X), \ n \in \N$.

\bigskip

\bigskip

We note that a $\SR$ flow $ (X,T) $ is $\D$, so its enveloping semigroup $E(X)$ is a group and the flow $(E(X),T)$ is minimal. We have more here, analogous to the result in \cite{TDES}:

\begin{theorem} \label{wr=Ett} For the flow $(X,T)$ the following are equivalent:
	
	\begin{enumerate}
		\item $(X,T)$ is weakly rigid ($\mathbf{WR}$).
		\item The identity $e $ is not isolated in $E(X)$, and the flow $(E(X), T)$  is topologically transitive.
	\end{enumerate}
	
\end{theorem}

\bo Consider a basic open set $B = [x_1, \ldots, x_n; U_1, \ldots, U_n]$ in the point-open topology on $E(X)$ such that $e \in B$. Since $(X,T)$ is weakly rigid, there exists $ t (\neq e)\in T$ such that $t(x_i) \in U_i $ i.e. $t \neq e \in B$ and so $e$ is not isolated in $E(X)$. Thus $E(X)$ is a perfect space and since the orbit ${Te}$ is dense in $E(X)$, the system $(E(X), T)$ is topologically transitive.

Note that since $e \in E(X)$ is not isolated, there exists $ t (\neq e)\in T$, such that $t \in B = [x_1, \ldots, x_n; U_1, \ldots, U_n]$, for any basic open set in $E(X)$. Thus $t(x_i) \in U_i $, and so $(X,T)$ is weakly rigid. \eo

\bigskip

By Proposition \ref{sr-wr} together with  Theorem \ref{wr=Ett}, we have 

\bc  $(X,T)$ is $\SR$ $\Longrightarrow$ $(E(2^X), T)$ is topologically transitive. \ec

\bigskip

We summarize the differences induced by  the classes of $\EQ$, $\SR$, and $\D$ flows in the table below:

$$ \begin{array}{|c|c|c|c|}
	\hline 
	
	(X,T) \text{\ is}  & \EQ  & \SR \setminus \EQ  & \D \setminus \SR \\
	\hline

	(2^X,T) \text{\ is}  & \text{{ equicontinuous }} & \text{{ weakly rigid }} & \text{{ NOT weakly rigid }} \\
	
	\hline
	
	(2^X,T) \text{\ is}  & \text{{ pointwise almost periodic }} & \text{{ densely almost periodic }} & \text{{ densely almost periodic  }} \\
	
	\hline 
	
	(2^X,T) \text{\ has}  & \text{{ every point recurrent}} & \text{{ every point recurrent}} & \text{{dense recurrent points}} \\
	
	\hline 
	
	(2^X,T) \text{\ has}  &  \Delta = \cR = RP & \Delta = \cR \subsetneqq RP & \emptyset  = \cR  \subsetneqq RP\\
	
	\hline 
	
	(\it{F}_{n}(X),T) \text{\ is}  & \EQ  & \SR \setminus \EQ  & \D \setminus \SR \\
	\hline

	(E(2^X),T) \text{\ is}  & \text{minimal} & \text{transitive but not minimal} & \text{NOT transitive} \\
	
	\hline 
	
\end{array}$$

\bigskip

\bigskip

\end{section}

\begin{section}{ Semiflows}
	
We recall the study on enveloping semigroups for a semiflow $(X,S)$ from \cite{TDES}.

Here we have an analogous definition for semiflows:
	
	\bd A semiflow (X, S) is said to be \emph{strongly rigid} if for every minimal idempotent $ u \in E(X, S) $  is an accumulation point in  $E(X, S)$ when endowed with  the compact-open topology, i.e. 
	there exists a (non-trivial) net $ \{s_k\} $ in $S$ such that the convergence $ s_k \to u $ 	is uniform. \ed
	
	We note that since every minimal idempotent $ u \in E(X, S) $  is an accumulation point in  $E(X, S)$ when endowed with  the compact-open topology, each minimal idempotent $ u \in E(X, S) $ is continuous.

	\br All equicontinuous semiflows are strongly rigid. 
Thus,
	
\centering{$\EQ \ \text{semiflows} \ \subsetneqq \ \SR \ \text{semiflows}$} \er
	
\br Note that strongly rigid semiflows need not be distal. An easy example is an aymptotic, equicontinuous semiflow. Thus, the identity  need not be a minimal idempotent in the enveloping semigroup of a strongly rigid semiflow. Here minimal idempotents could also be constant mappings (say for an equicontinuous, proximal semiflow). \er	

\br Note that a finite product of $\SR$ semiflows will also be $\SR$ and a factor of $\SR$ semiflow will be $\SR$. The proof will be exactly on lines of that for flows. \er

We skip any discussion in general on strongly rigid semiflows and concentrate on the case of distal semiflows.

\subsection{Distal Semiflows}	We recall from \cite{TDES}:

If each $s \in S$ is surjective,  then the semiflow $(X, S)$ is called a \emph{surjective semiflow}. If each $s \in S$ is bijective, then \emph{$[S]$} denotes the smallest group of self-homeomorphisms of $X$ containing $S$. In such a case the semiflow $(X, S)$ induces the flow $(X,[S])$. A point $x \in X$ is called a \emph{distal point} of $(X,S)$ if $x$ is the only point proximal to itself in $\overline{Sx}$. If each $x \in X$ is distal then the semiflow $(X,S)$ is called distal.
	
\begin{proposition} \cite{TDES}
	If $(X,S)$ is a distal semiflow, then $(X,S)$ is surjective.
\end{proposition}	
	
Distal systems are always  injective, and hence we note that distal systems are always bijective. Thus every distal semiflow $(X,S)$ induces a distal flow $(X,[S])$.

And, we have

\begin{theorem} \cite{TDES}
	For a semiflow $(X,S)$ the following conditions are equivalent:
	\begin{enumerate}
		\item $(X,S)$ is distal.
		\item $(E(X,S), S)$ is minimal.
		\item $E(X,S)$ is a group with identity $e$ being the only idempotent.
		
		Furthermore, $E(X,S) = E(X,[S])$.
	\end{enumerate}
	
\end{theorem}

\bex

Recall Example \ref{A1}. And for $X$ and $f$ considered there, consider the semicascade $(X,f)$. This semicascade is distal but not strongly rigid.

\medskip

Recall Example \ref{A2}. For the map $T: X \to X$ considered there, let $S =\{T^s: s=3,4,5, \ldots\}$, then $S$ is a semigroup under the operation of composition of mappings on $X$. Thus $(X,S)$ gives a semiflow. Note that $(X,S)$ is distal and so induces the distal flow 
$(X,[S])$ which is same as the cascade $(X,T)$ considered there.

Recalling the details discussed in Example \ref{A2}, we can say the the semiflow $(X,S)$ is strongly rigid, and not equicontinuous.
\eex	
	
Thus, in the category of distal semiflows we have $  \SR \ \text{semiflows} \ \subsetneqq \ \D \ \text{semiflows}$.

\bp A distal $\SR$ semiflow $(X,S)$ generates a $\SR$ flow $(X,{[S]})$. \ep	
	
Consider a distal semiflow $(X,S)$. Let $[S]$ be the group generated by the semigroup $S$.   Note that for a net $\{s_\alpha\}$ in $[S]$ such that $s_\alpha \longrightarrow e$, there exists a subnet $\{s_{\alpha_k}\}$ in either $S$ or $S^{-1}$ such that $s_{\alpha_k} \longrightarrow e$. Thus, we can say

\bp For a semigroup $S$ generating the group $[S]$, 

 The $\SR$ flow $(X,[S])$ generates either the $\SR$ semiflow $(X,S)$ or $(X,S^{-1})$. 
\ep

Let $[S]	= S \cup S^{-1}$.  Since $[S]$  is a topological group, $s_{\alpha} \longrightarrow e \Rightarrow s_{\alpha}^{-1} \longrightarrow e$. 

\bp If $[S]	= S \cup S^{-1}$, the distal semiflow $(X,S)$ is $\SR$ $\Leftrightarrow$  the distal semiflow $(X,S^{-1})$ is $\SR$.\ep
	
	\end{section}

\vspace{12pt}
\bibliography{xbib}

\begin{thebibliography}{99}
	
	\bibitem {AAA}
	{\bf Ethan Akin, Joseph Auslander and Anima Nagar,} Dynamics of Induced Systems, {\it Ergod. Th. and Dynam. Sys.,} 37 (2017), 2034-2059.
	
	\bibitem {AUS}
	{\bf Joseph Auslander,} Minimal flows and their extensions, {\it North-Holland Mathematics studies,} 153 (1988).
	
	\bibitem {ELL}
	{\bf David B. Ellis and Robert Ellis,} Automorphisms and Equivalence Relations in Topological Dynamics, {\it Cambridge University Press,} (2014).
	
	
	\bibitem {RIG}
	{\bf Eli Glasner and David Maon,} Rigidity in topological dynamics, {\it Ergod. Th. and Dynam. Sys.,} 9(1989), 309-320.
	
	
	\bibitem{GH} 
	{\bf W. H. Gottschalk and G. A. Hedlund}, Topological dynamics, {\it Amer.
	Math. Soc. Colloquium Publications}, vol. 36, (1955).


\bibitem{OP}
{\bf Jie Li, Piotr Oprocha, Xiangdong Ye, Ruifeng Zhang,} When are all closed subsets recurrent?, {\it Ergod. Th. and Dynam. Sys.,}  37(2017), 2223-2254.



	
	\bibitem {TDES}
	{\bf Anima Nagar and Manpreet Singh,} Topological Dynamics of the Enveloping Semigroup, {\it  arXiv:1810.12854v2} (2019).
	
	
	\bibitem {AN}
	{\bf Anima Nagar,} Characterization of Quasifactors, {\it arXiv:2006.11611} (2020).
	
	
\end{thebibliography}

\end{document}